\newcount\secno
\newcount\prmno
\newif\ifnotfound
\newif\iffound

\def\namedef#1{\expandafter\def\csname #1\endcsname}
\def\nameuse#1{\csname #1\endcsname}

\long\def\ifundefined#1#2#3{\expandafter\ifx\csname
  #1\endcsname\relax#2\else#3\fi}
\def\hwrite#1#2{{\let\the=0\edef\next{\write#1{#2}}\next}}

\toksdef\ta=0 \toksdef\tb=2
\long\def\leftappenditem#1\to#2{\ta={\\{#1}}\tb=\expandafter{#2}%
                                \edef#2{\the\ta\the\tb}}
\long\def\rightappenditem#1\to#2{\ta={\\{#1}}\tb=\expandafter{#2}%
                                \edef#2{\the\tb\the\ta}}

\def\lop#1\to#2{\expandafter\lopoff#1\lopoff#1#2}
\long\def\lopoff\\#1#2\lopoff#3#4{\def#4{#1}\def#3{#2}}

\def\ismember#1\of#2{\foundfalse{\let\given=#1%
    \def\\##1{\def\next{##1}%
    \ifx\next\given{\global\foundtrue}\fi}#2}}

\def\section#1{\vskip1cm
               \global\def\currenvir{section}
               \global\advance\secno by1\global\prmno=0
               {\bf \number\secno. {#1}}
               \smallskip}

\def\subsection{\global\def\currenvir{subsection}
                \global\advance\prmno by1
                \smallskip\ind{(\number\secno.\number\prmno) }}
\def\subsec{\global\def\currenvir{subsection}
                \global\advance\prmno by1
                { (\number\secno.\number\prmno)\ }}

\def\proclaim#1{\global\advance\prmno by 1
                {\bf #1 \the\secno.\the\prmno$.-$ }}

\long\def\th#1 \enonce#2\endth{%
   \medbreak\proclaim{#1}{\it #2}\global\def\currenvir{th}\smallskip}

\def\rem#1{\global\advance\prmno by 1
\smallskip {\it #1 }\  \the\secno.\the\prmno$.-$}


\def\isinlabellist#1\of#2{\notfoundtrue%
   {\def\given{#1}%
    \def\\##1{\def\next{##1}%
    \lop\next\to\za\lop\next\to\zb%
    \ifx\za\given{\zb\global\notfoundfalse}\fi}#2}%
    \ifnotfound{\immediate\write16%
                 {Warning - [Page \the\pageno] {#1} No reference found}}%
                \fi}%
\def\ref#1{\ifx\labellist\empty{\immediate\write16
                 {Warning - No references found at all.}}
               \else{\isinlabellist{#1}\of\labellist}\fi}

\def\newlabel#1#2{\rightappenditem{\\{#1}\\{#2}}\to\labellist}
\def\labellist{}
\def\label#1{%
  \def\given{th}%
  \ifx\given\currenvir%
    {\hwrite\lbl{\string\newlabel{#1}{\number\secno.\number\prmno}}}\fi%
  \def\given{section}%
  \ifx\given\currenvir%
    {\hwrite\lbl{\string\newlabel{#1}{\number\secno}}}\fi%
  \def\given{subsection}%
  \ifx\given\currenvir%
    {\hwrite\lbl{\string\newlabel{#1}{\number\secno.\number\prmno}}}\fi%
  \def\given{subsubsection}%
  \ifx\given\currenvir%
  {\hwrite\lbl{\string%
    \newlabel{#1}{\number\secno.\number\subsecno.\number\subsubsecno}}}\fi
  \ignorespaces}
\newwrite\lbl
\def\openall{\openout\lbl=\jobname.lbl}

\newread\testfile
\def\lookatfile#1{\openin\testfile=\jobname.#1
    \ifeof\testfile{\immediate\openout\nameuse{#1}\jobname.#1
                    \write\nameuse{#1}{}
                    \immediate\closeout\nameuse{#1}}\fi%
    \immediate\closein\testfile}%

\def\begin{\newlabel{def}{1.1}
\newlabel{qt}{2.3}
\newlabel{nil}{2.5}
\newlabel{loc}{3.2}
\newlabel{W}{3.3}}


\magnification 1250
\pretolerance=500 \tolerance=1000  \brokenpenalty=5000
\mathcode`A="7041 \mathcode`B="7042 \mathcode`C="7043
\mathcode`D="7044 \mathcode`E="7045 \mathcode`F="7046
\mathcode`G="7047 \mathcode`H="7048 \mathcode`I="7049
\mathcode`J="704A \mathcode`K="704B \mathcode`L="704C
\mathcode`M="704D \mathcode`N="704E \mathcode`O="704F
\mathcode`P="7050 \mathcode`Q="7051 \mathcode`R="7052
\mathcode`S="7053 \mathcode`T="7054 \mathcode`U="7055
\mathcode`V="7056 \mathcode`W="7057 \mathcode`X="7058
\mathcode`Y="7059 \mathcode`Z="705A
\def\spacedmath#1{\def\packedmath##1${\bgroup\mathsurround =0pt##1\egroup$}
\mathsurround#1
\everymath={\packedmath}\everydisplay={\mathsurround=0pt}}
\def\nospacedmath{\mathsurround=0pt
\everymath={}\everydisplay={} } \spacedmath{2pt}
\def\qfl#1{\buildrel {#1}\over {\longrightarrow}}

\def\hfl#1#2{\normalbaselines{\baselineskip=0truept
\lineskip=10truept\lineskiplimit=1truept}\nospacedmath\smash{\mathop{\hbox to
12truemm{\rightarrowfill}}\limits^{\scriptstyle#1}_{\scriptstyle#2}}}
\def\diagram#1{\def\normalbaselines{\baselineskip=0truept
\lineskip=10truept\lineskiplimit=1truept}   \matrix{#1}}
\def\vfl#1#2{\llap{$\scriptstyle#1$}\left\downarrow\vbox to
6truemm{}\right.\rlap{$\scriptstyle#2$}}

\def\mono{\lhook\joinrel\mathrel{\longrightarrow}}
\def\iso{\mathrel{\mathop{\kern 0pt\longrightarrow }\limits^{\sim}}}

\def\union_#1{\raise
2pt \hbox{$\mathrel{\mathop{\kern0pt{\scriptscriptstyle\bigcup}}
\limits_{#1}}$}}
\def\sdir_#1^#2{\mathrel{\mathop{\kern0pt\oplus}\limits_{#1}^{#2}}}
\def\pprod_#1^#2{\raise
2pt \hbox{$\mathrel{\scriptstyle\mathop{\kern0pt\prod}\limits_{#1}^{#2}}$}}
\font\eightrm=cmr8         \font\eighti=cmmi8
\font\eightsy=cmsy8        \font\eightbf=cmbx8
\font\eighttt=cmtt8        \font\eightit=cmti8
\font\eightsl=cmsl8        \font\sixrm=cmr6
\font\sixi=cmmi6           \font\sixsy=cmsy6
\font\sixbf=cmbx6\catcode`\@=11
\def\eightpoint{%
  \textfont0=\eightrm \scriptfont0=\sixrm \scriptscriptfont0=\fiverm
\def\rm{\fam\z@\eightrm}%
  \textfont1=\eighti  \scriptfont1=\sixi 
\scriptscriptfont1=\fivei
  \textfont2=\eightsy \scriptfont2=\sixsy \scriptscriptfont2=\fivesy
  \textfont\itfam=\eightit
  \def\it{\fam\itfam\eightit}%
  \textfont\slfam=\eightsl
  \def\sl{\fam\slfam\eightsl}%
  \textfont\bffam=\eightbf \scriptfont\bffam=\sixbf
  \scriptscriptfont\bffam=\fivebf
  \def\bf{\fam\bffam\eightbf}%
  \textfont\ttfam=\eighttt
  \def\tt{\fam\ttfam\eighttt}%
  \abovedisplayskip=9pt plus 3pt minus 9pt
  \belowdisplayskip=\abovedisplayskip
  \abovedisplayshortskip=0pt plus 3pt
  \belowdisplayshortskip=3pt plus 3pt 
  \smallskipamount=2pt plus 1pt minus 1pt
  \medskipamount=4pt plus 2pt minus 1pt
  \bigskipamount=9pt plus 3pt minus 3pt
  \normalbaselineskip=9pt
  \normalbaselines\rm}\catcode`\@=12
\newcount\noteno
\noteno=0
\def\up#1{\raise 1ex\hbox{\sevenrm#1}}
\def\note#1{\global\advance\noteno by1
\footnote{\parindent0.4cm\up{\number\noteno}\
}{\vtop{\eightpoint\baselineskip12pt\hsize15.5truecm\noindent
#1}}\parindent 0cm}

\font\san=cmssdc10
\def\ext{\hbox{\san \char3}}

\def\pc#1{\tenrm#1\sevenrm}
\def\tx{\kern-1.5pt -}
\def\cqfd{\kern 2truemm\unskip\penalty 500\vrule height 4pt depth 0pt width
4pt\medbreak} 
\def\virg{\raise
.4ex\hbox{,}}
\def\decale#1{\smallbreak\hskip 28pt\llap{#1}\kern 5pt}
\def\no{n\up{o}\kern 2pt}
\def\ind{\par\hskip 1truecm\relax}
\def\indp{\par\hskip 0.5truecm\relax}
\def\moins{\mathrel{\hbox{\vrule height 3pt depth -2pt width 6pt}}}
\def\rond{\kern 1pt{\scriptstyle\circ}\kern 1pt}
\def\iso{\mathrel{\mathop{\kern 0pt\longrightarrow }\limits^{\sim}}}

\def\Hom{\mathop{\rm Hom}\nolimits}

\def\det{\mathop{\rm det}\nolimits}

\def\div{\mathop{\rm div\,}\nolimits}

\def\codim{\mathop{\rm codim}\nolimits}
\def\Sing{\mathop{\rm Sing}\nolimits}
\input amssym.def
\input amssym
\vsize = 25truecm
\hsize = 16truecm
\voffset = -.5truecm
\parindent=0cm
\baselineskip15pt
\overfullrule=0pt
\begin
\vskip1cm
\centerline{\bf  Symplectic singularities}
\smallskip
\centerline{Arnaud {\pc BEAUVILLE}} 
\vskip1cm
{\bf Introduction}
\smallskip
\ind We introduce in this paper  a particular class of rational
singularities, which we call {\it symplectic}, and classify the simplest
ones. Our motivation comes from the  analogy between rational
 Gorenstein singularities and Calabi-Yau manifolds: a compact, K\"ahler
manifold of dimension $n$ is a Calabi-Yau manifold
if it admits a nowhere vanishing
$n$\tx form, while  a normal  variety $V$ of dimension
$n$ has rational Gorenstein singularities\note{also called canonical
singularities of index $1$.} if its smooth part
$V_{\rm reg}$ carries a nowhere vanishing
$n$\tx form, with the extra property that its pull-back in any resolution
$X\rightarrow V$ extends to a holomorphic form on $X$. Among Calabi-Yau
manifolds an important role is played by the  symplectic (or hyperk\"ahler)
manifolds, which admit a
 holomorphic, everywhere non-degenerate $2$\tx form; by analogy we say that a 
normal  variety $V$ has {\it symplectic singularities}  if $V_{\rm reg}$
carries a closed symplectic $2$\tx form whose pull-back in any resolution
$X\rightarrow V$ extends to a holomorphic $2$\tx form on $X$.
\ind  We will look for the simplest possible isolated symplectic
singularities ${\rm o}\in V$, namely those whose projective tangent cone is
smooth: this means that  blowing up  ${\rm o}$ in $V$ provides a resolution of
$V$ with a smooth exceptional divisor. Examples of such
singularities  are obtained as
follows. Each  simple complex Lie algebra has a smallest non-zero nilpotent
orbit 
${\cal O}_{\rm min}$  for the adjoint action; its closure $\overline{\cal
O}_{\rm min}={\cal O}_{\rm min}\cup\{0\}$ has a symplectic singularity at $0$,
isomorphic to the cone over the smooth variety ${\bf P}{\cal
O}_{\rm min}:={\cal O}_{\rm min}/{\bf C^*}$. In
particular its  projective tangent cone is smooth (it is isomorphic to 
${\bf P}{\cal O}_{\rm min}$).
\ind Our main result is the converse:\smallskip 
{\bf Theorem}$.-$
{\it  Let $(V,{\rm o})$ be an isolated symplectic singularity, whose
projective tangent cone is smooth. Then $(V,{\rm o})$ is analytically
isomorphic to
$(\overline{\cal O}_{\rm min},0)$ for some simple complex Lie algebra}.
\smallskip 
\ind The key point of the proof is the fact that the homogeneous space 
${\bf P}{\cal O}_{\rm min}$ carries a  holomorphic {\it contact structure}
(inherited from the symplectic structure of ${\cal O}_{\rm min}$). Given a
resolution $X\rightarrow V$ with a smooth exceptional divisor $E$, we show 
that the extension to $X$ of the symplectic form has a
residue on $E$ which defines a contact structure. We then deduce from [B1]
 that $E$ is isomorphic to some ${\bf P}{\cal O}_{\rm min}$,
and we conclude with a classical criterion of Grauert. 
\ind We discuss in \S 4 whether a classification of isolated symplectic
singularities makes sense. Each such singularity gives rise to many others by
considering its quotient by a finite group; to get rid of those we propose to
consider only isolated symplectic singularities with trivial local fundamental
group. The singularities $(\overline{\cal O}_{\rm min},0)$ have this property
when the Lie algebra is not of type $C_l$; it is certainly desirable to find
more examples.

\section{Definition and basic properties}
\ind We consider algebraic varieties over ${\bf C}$ (our results extend
readily to the analytic category). We
will say  that a holomorphic 2-form on a smooth variety is {\it
symplectic} if it is closed and non-degenerate at every point. A {\it
resolution} of an algebraic variety $V$ is a proper, birational morphism
$f:X\rightarrow V$ where $X$ is smooth.
\th Definition
\enonce A  variety has a symplectic singularity at a point if this point 
admits an open neighborhood $V$ such that:
\indp {\rm a)} $V$ is normal;
\indp {\rm b)} The smooth part $V_{\rm reg}$ of $V$ admits a symplectic $2$\tx
form
$\varphi$;
\indp {\rm c)} For any resolution $f:X\rightarrow V$, the pull back of
$\varphi$ to $f^{-1} (V_{\rm reg})$ extends to a holomorphic $2$\tx form on
$X$.
\endth\label{def}
\ind We will mostly consider a symplectic singularity as a germ $(V,{\rm o})$
-- in which case we will always assume that $V$ satisfies the conditions a)
to c).
\subsection As for rational singularities it is enough to check c) for one
particular resolution: the point is that if $X\rightarrow Y$ is a proper
birational morphism of smooth varieties, a meromorphic form $\varphi$ on $Y$
is holomorphic if and only if $g^*\varphi$ is holomorphic.
\th Proposition
\enonce A symplectic singularity is rational Gorenstein.
\endth
\indp{\it Proof}: We keep the notation of Definition \ref{def} and put $\dim
V=2r$. The form $\varphi^r$ generates the line bundle
$\omega^{}_{V_{\rm reg}}$, and for any resolution $X\rightarrow V$ extends to
a holomorphic form on $X$; this implies that $V$ has  rational Gorenstein
singularities [R].\cqfd
\ind The following remark shows that isolated symplectic singularities of 
dimension $>2$ are {\it not} local complete intersections:
\th Proposition
\enonce Let $V$ be a variety with symplectic singularities which is locally 
a complete intersection. Then the singular locus of $V$ has codimension
$\le 3$.
\endth
\indp{\it Proof}:  We can realize
locally
$V$ as a complete intersection in some smooth variety $S$. 
The exact sequence
$$0\rightarrow N_{V/S}^* \longrightarrow
\Omega^1_S{}^{}_{|V}\longrightarrow
\Omega^1_V\rightarrow 0\ ,$$
provides a length $1$ locally free resolution of $\Omega^1_V$. We can assume
$\codim \Sing(V)\ge 3$; by the Auslander-Buchsbaum theorem and the fact that
$V$ is Cohen-Macaulay,  the depth of $\Omega^1_V$ at every point of
$\Sing(V)$ is $\ge 2$.  It follows that $\Omega^1_V$ is a reflexive sheaf,
so the isomorphism
$\Omega^1_{V_{\rm reg}}\rightarrow T^{}_{V_{\rm reg}}$ defined by a symplectic
$2$\tx form on $V_{\rm reg}$
 extends to an isomorphism $\Omega^1_V\rightarrow T^{}_V$.
Combining the resolution of $\Omega^1_V$ and its dual we get an exact
sequence
$$0\rightarrow N_{V/S}^* \longrightarrow
\Omega^1_S{}^{}_{|V}\longrightarrow  T_S{}^{}_{|V}\qfl{u} N_{V/S}\ ,$$
 where the support of the cokernel $T^1$ of $u$ is exactly $\Sing(V)$.
Using  the Auslander-Buchsbaum theorem again we get
$\dim(T^1)=\dim\Sing (V)\ge
\dim(V)-3$.\cqfd
 
\section{Examples}
\subsection In dimension $2$, the symplectic singularities are the rational
double points (that is the A-D-E singularities).\par
\subsection Any product of varieties with symplectic singularities has again
symplectic singularities. 
\medskip
\subsec {\it Quotient singularities}\label{qt}
\ind The following result will provide us with a large list of symplectic
singularities:
\th Proposition
\enonce Let $V$ be a variety with symplectic singularities, 
$G$ a finite group of automorphisms of $V$, preserving a symplectic $2$\tx
form on $V_{\rm reg}$.  Then $V/G$ has symplectic singularities.
\endth
\indp{\it Proof}: We first observe that the fixed locus  $F_g$ in $V_{\rm
reg}$ of any element $g\not=1$ in $G$ is a symplectic subvariety of
$V_{\rm reg}$, and therefore has codimension $\ge 2$. Let
$ V^{\rm o}:=V_{\rm reg}\moins \union_{g\not=1}F_g$. The symplectic
2-form on $V^{\rm o}$ descends to a symplectic 2-form  on
$V^{\rm o}/G$, which extends to a symplectic 2-form $\varphi$ on
$(V/G)_{\rm reg}$. Let $g:Y\rightarrow V/G$ be a resolution of $V/G$; we
can find a commutative diagram
$$\diagram{X & \hfl{f}{} & V \cr
\vfl{}{} & & \vfl{}{} \cr
Y & \hfl{g}{} & V/G \cr
}$$where $f$ is a resolution of $V$. 
Then  $g^*\varphi$ is a meromorphic 2-form on $Y$, whose pull back to
$X$ is holomorphic. By an easy local computation,  this implies that
$g^*\varphi$ is holomorphic.\cqfd
 \medskip
\ind This applies for instance when  $V$ is a finite-dimensional
symplectic vector space, and  $G$ a finite subgroup of ${\rm Sp}(V)$.
If we impose moreover that the non trivial elements of $G$ have all their
eigenvalues $\not=1$, then $V/G$ has an isolated symplectic
singularity. The case $\dim(V)=2$ gives the rational double points.  
\ind The Proposition also applies to the symmetric products
$V^{(p)}=V^p/{\goth S}_p$: if $V$ has symplectic singularities, so does
$V^{(p)}$. 
\medskip 
\subsec {\it Nilpotent orbits}\label{nil}
\ind  Let ${\goth g}$ be a simple complex Lie
algebra
 and ${\cal O}\i{\goth g}$ a  nilpotent orbit (for the
adjoint action)\note{A general reference for  nilpotent
orbits is [C-M].}. Then {\it the normalization of the closure
 of ${\cal O}$ in ${\goth g}$ has symplectic singularities}. This is due to
Panyushev [P], who uses it to prove that this variety has rational Gorenstein
singularities. The point is that ${\cal O}$ can be identified with a
coadjoint orbit using the Killing form, and therefore carries the
Kostant-Kirillov symplectic
$2$\tx form.
\ind In particular, the Lie algebra ${\goth g}$ contains a unique
(non-zero) minimal nilpotent orbit  ${\cal O}_{\rm min}$, which is contained
in the closure of all non-zero nilpotent orbits. The closure 
$\overline{\cal O}_{\rm min}= {\cal O}_{\rm min}\cup\{0\}$ is normal, and has
an isolated symplectic singularity at
$0$. 
\ind This singularity can be described as follows. The orbit
${\cal O}_{\rm min}$ is stable by homotheties; the quotient 
${\bf P}{\cal O}_{\rm min}:={\cal O}_{\rm min}/{\bf C}^*$ is a smooth,
closed subvariety of ${\bf P}({\goth g})$.
 The variety $\overline{\cal O}_{\rm min}$ is the cone over ${\bf P}{\cal
O}_{\rm min}\i {\bf P}({\goth g})$. This means that we have a resolution
$f:L^{-1}
\rightarrow \overline{\cal O}_{\rm min}$, where $L$ is the restriction
 of  ${\cal O}_{{\bf P}({\goth g})}(1)$ to ${\bf P}{\cal O}_{\rm min}$,
and $f$  contracts to $0$ the  zero section $E$ of $L^{-1} $. In this
situation $f$ is the blow up of
$0$ in $\overline{\cal O}_{\rm min}$, and the exceptional divisor $E$,
isomorphic to ${\bf P}{\cal O}_{\rm min}$, is the projective tangent cone
to $0$ in
$\overline{\cal O}_{\rm min}$. 
 
\section{Characterization of minimal orbits singularities}
\subsection This section is devoted to the proof of the theorem stated in the
introduction. So we let $(V, {\rm o})$ be an isolated symplectic singularity,
  $f:X\rightarrow V$  the blow up of the
maximal ideal of ${\rm o}$ in $V$, and $E$ the exceptional
divisor. By construction $E$ is isomorphic to the projective tangent cone
to $V$ at ${\rm o}$; we assume that it  is smooth. Since
$E$ is a Cartier divisor in $X$ it follows that $X$ is
smooth. 
\ind We denote by $i$  the embedding
of $E$ in $X$, and put $L:=i^*{\cal O}_X(-E)$.
 By the standard properties of the blow up the line bundle $L$ on $E$ is
{\it very ample}.  
 \subsection \label{loc} Let $\dim V=2r$. We can assume that $V\moins
\{{\rm o}\}$ carries a symplectic
$2$\tx form  which extends to a holomorphic $2$\tx form $\varphi$ on $X$;
we have $\div(\varphi^r)=kE$ for some integer $k\ge 0$. The adjunction
formula gives $K_E=L^{-k-1} $, so that $E$ is a Fano manifold. This implies 
$H^0(E,\Omega^p_E)=0$ for each $p\ge 1$, and in particular $i^*\varphi=0$.
\ind Let $e\in E$. Since $\varphi$ is closed,  we can write $\varphi=d\alpha$
 in a neighbourhood $U$ of $e$ in $X$, where $\alpha$ is a $1$\tx form on
$U$ such that $i^*\alpha$ is closed. Shrinking $U$ if necessary we can
write $i^*\alpha=d(i^*g)$ for some function $g$ on $U$; replacing $\alpha$
by $\alpha-dg$ we may assume $i^*\alpha=0$. If $u=0$ is a local equation
of $E$ in $U$, this means that $\alpha$ is of the form $u\,\tilde
\theta+h\,du$, where $\tilde\theta$ is a $1$\tx form and $h$ a function on
$U$; replacing $\alpha$ by $\alpha-d(hu)$ and $\tilde \theta$ by $\tilde
\theta-dh$ we arrive at
$\alpha=u\tilde\theta$ and
$$\varphi=du\wedge\tilde\theta + u\,d\tilde\theta\ .$$

\ind This gives $\varphi^r=ru^{r-1}du\wedge
\tilde \theta\wedge (d\tilde \theta)^{r-1}+u^r(d\tilde
\theta)^r$. Thus the order of vanishing
$k$ of $\varphi^r$ along $E$ is 
$\ge r-1$; the crucial point of the proof is the equality
$k=r-1$. We need an easy lemma:

\th Lemma
\enonce  Let $X$ be a smooth closed
submanifold of a projective space ${\bf P}^N$, of  degree
$\ge 2$. Then
$H^0(X,\ext^pT_X(-p))=0$ for
$0<p<\dim(X)$, and for $p=\dim(X)$ except if $X$ is a
hyperquadric.
\endth\label{W}

\indp{\it Proof}: When $X$ is a hyperquadric our assertion is equivalent to 
$H^0(X,\Omega^q_X(q))=0$ for $0<q<\dim(X)$, which can be checked by a
direct computation (see for instance [K], thm. 3). We assume
$\deg(X)\ge 3$.
\ind The case $p=1$ follows from a more
general result of Wahl ([W], see remark below). 
Then we use induction on the dimension of
$X$, the case of curves being clear. Let
$H$ be a smooth hyperplane section of $X$; the exact sequence
$$0\rightarrow T_H\longrightarrow T_X{}^{}_{|H}\longrightarrow
{\cal O}_H(1)\rightarrow 0$$
gives rise to exact sequences
$$0\rightarrow \ext^pT_H\,(-p)\longrightarrow
\ext^pT_X{}^{}_{|H}(-p)\longrightarrow \ext^{p-1}T_H(-(p-1))\rightarrow 0$$
By the induction hypothesis we conclude that $H^0(X,\ext^pT_X{}^{}_{|H}(-p))$
is zero for $ p\ge 2$, and therefore $H^0(X,\ext^pT_X{}^{}(-p))=0$
 for $p\ge 2$.\cqfd

\rem{Remark} Wahl's result is rather
easy in our situation: using the exact sequence 
$$0\rightarrow H^0(X,T_X(-1))\longrightarrow  H^0(X,T_{{\bf
P}^N}(-1)^{}_{|X})\longrightarrow H^0(X,N_{X/{\bf P}^N}(-1))$$
 and the
isomorphism ${\bf C}^{N+1}\iso H^0(X,T_{{\bf
P}^N}(-1)^{}_{|X})$ deduced from the Euler exact sequence, we see that a
nonzero element of $H^0(X,T_X(-1))$ corresponds to a point $p\in{\bf
P}^N$ such that all projective  tangent spaces
${\bf P}T_x(X)$, for  $x$ in $X$,  pass through $p$. This is easily seen
to be impossible, for instance by induction on $\dim(X)$. 
\ind It seems natural to conjecture that the statement of the lemma extends
to the more general situation considered in [W], namely that 
$H^0(X,\ext^pT_X\otimes L^{-p})=0$ {\it  for $p>0$ whenever $L$ is ample,
except if $(X,L)=({\bf P}^n,{\cal O}_{{\bf P}^n}(1))$, with $n\ge p$, or} 
$(X,L)=(Q_p, {\cal O}_{Q_p}(1))$.
 \smallskip 
\subsection We now prove the equality $k=r-1$. If
$E={\bf P}^{2r-1}$ and $L={\cal O}_{{\bf P}^{2r-1}}(1)$,  $V$ is
smooth; if $E={\bf P}^{1}$ and $L={\cal O}_{{\bf P}^{1}}(2)$, $V$ is a
surface with an ordinary double point. We exclude these two cases. The perfect
pairing 
$\Omega^1_X\otimes\Omega^{2r-1}_X\rightarrow K_X$ provides an isomorphism
$\Omega^{2r-1}_X\cong T_X\otimes K_X$; thus 
exterior product with
$\varphi^{r-1}$ gives a linear map 
$\Omega^1_X\rightarrow T_X(kE)$, which is an isomorphism
outside $E$ (it is the inverse of the isomorphism defined by $\varphi$).
This map may vanish on
$E$, say with order
$k-j$
$(j\le k)$, so that we get a map $\lambda:\Omega^1_X\rightarrow T_X(jE)$
whose restriction to $E$ is nonzero. Observe that $\det \lambda$ is a
section of 
${\cal O}_X(2(rj-k)E)$ which is nonzero outside $E$, hence $k\le rj$ and in 
particular $j\ge 0$. 
\def\lr{\ \longrightarrow \ }
\ind We have a diagram of exact sequences
$$\diagram{& 0\lr  L \lr 
\Omega^1_X{}^{}_{|E}
\lr \Omega^1_E\lr 0 &\cr
&\vfl{}{\lambda_{|E}}&\cr
&0\lr  T _E\otimes
L^{-j}  \lr T_X{}^{}_{|E}\otimes L^{-j}\lr 
L^{-j-1}\lr 0 &\ .}$$
Since $j\ge 0$ we have $\Hom(L,L^{-j-1})=\Hom(\Omega^1_E,L^{-j-1})=\Hom(L,T
_E\otimes L^{-j})=0$ by lemma \ref{W}. Thus $\lambda^{}_{|E}$ factors through 
a map $\mu: \Omega^1_E\rightarrow T _E\otimes L^{-j}$; since $\lambda$ is
antisymmetric $\mu$ is antisymmetric, that is comes from an element of 
$H^0(E,\ext^{2}T_E\otimes L^{-j})$. 
\ind Since $\lambda^{}_{|E}$ is non-zero, lemma \ref{W} implies $j\le 1$,
hence
$k\le r$. Moreover if $k=rj$,  $\det\lambda$ does not vanish, hence $\lambda$ 
and therefore $\lambda^{}_{|E}$ are isomorphisms; but this is impossible
because
$\lambda^{}_{|E}$ vanishes on the sub-bundle $L\i\Omega^1_X{}^{}_{|E}$. Thus
we have $k<rj$, and therefore $j=1$ and $k=r-1$.
\subsection Going back to the local computation of (\ref{loc}),  we observe
that the form $\theta:=i^*\tilde\theta$ is defined globally as a section
of $\Omega^1_E\otimes L$: it is the image of $\varphi\in
H^0(X,\Omega^2_X(\log E)(-E))$ by the residue map  $\Omega^2_X(\log
E)(-E)\rightarrow \Omega^1_E\otimes {\cal O}_X(-E)^{}_{|E}$. 
We now know that the $(2r)$\tx form
$du\wedge\tilde \theta\wedge (d\tilde \theta)^{r-1}$ on $U$ does not vanish,
so the twisted $(2r-1)$\tx form $\theta\wedge (d\theta)^{r-1}\in
H^0(E,K_E\otimes L^r)$ does not vanish. This means, by definition, that
$\theta$ is a {\it contact structure} on the Fano manifold $E$. The
classification of Fano contact manifolds is an interesting problem, with
important applications to Riemannian geometry (see for instance [L] or [B2]).
Here we have one more information, namely that the line bundle $L$ is {\it
very} ample; this implies that $E$ is isomorphic to one of the homogeneous
contact manifolds
${\bf P}{\cal O}_{\rm min}$ ([B1], cor. 1.8).
\subsection It remains to show that the embedding of $E$  in $X$ 
is isomorphic, in some open
neighbourhood of $E$, to the embedding of the zero section in the line bundle
$L^{-1}
\rightarrow E$. By a criterion of Grauert [G], it is sufficient to prove that
the spaces
$H^1(E,T_E\otimes L^k)$ and $H^1(E,L^k)$ are zero for $k\ge 1$. The second
assertion follows from the Kodaira vanishing theorem, since $E$ is a Fano
manifold. The first one   can be deduced (with some work) from the Bott
vanishing theorem, but we found more convenient to use the following easy
lemma:
\th Lemma
\enonce Let $Z$ be a compact complex manifold, $E$ a vector bundle on $X$
spanned by its global sections, $L$ an ample line bundle on $X$. Then
$H^p(Z,K_Z\otimes E\otimes\det E\otimes L)=0$ for $p\ge 1$.
\endth
\indp{\it Proof}: Put $P={\bf P}_Z(E)$, and let  $p:P\rightarrow Z$ be the
canonical fibration. Choose a surjective map  ${\cal
O}_Z^{N+1}\rightarrow E$. It induces an embedding $P\mono Z\times {\bf
P}^N$. For each $j\ge 1$ the line bundle $p^*L\otimes {\cal O}_P(j)$ is the
restriction  to $P$ of $L\boxtimes {\cal O}_{{\bf P}^N}(j)$, and therefore
is ample. We have $K_P\cong p^*(K_Z\otimes\det E)\otimes {\cal
O}_P(-r)$, with $r={\rm rk}(E)$. By the Kodaira vanishing theorem 
we have $H^p(P, p^*(K_Z\otimes\det E\otimes L)\otimes {\cal O}_P(1))=0$
for $p\ge 1$. Since $Rp_*(p^*M\otimes {\cal O}_P(1)))\cong M\otimes E$ for
any line bundle $M$ on $Z$, our assertion follows.\cqfd 

\section{Local fundamental group}
\subsection In view of (\ref{qt}) it seems hopeless to classify all isolated
symplectic singularities: there are too many quotient singularities, already
in dimension $4$. One way to get around this problem is to consider only
singularities with {\it trivial local fundamental group}. We briefly recall
the definition: if
$(V,{\rm o})$ is an isolated singularity, we can find a fundamental system
$(V_n)_{n\ge 1}$ of neighbourhoods of
${\rm o}$ such that $V_q$ is a deformation retract of $V_p$ for $q\ge p$;
the group $\pi _1(V_n)$, which is independant of $n$ and of the particular
fundamental system, is called   the local fundamental group of $V$ at ${\rm
o}$ and denoted $\pi _1^{\rm o}(V)$ (for a canonical definition one should be
more careful about base points, but this is irrelevant here).  
\ind If $(V,{\rm o})$ is a quotient of an isolated singularity $(W,\omega)$ by 
a finite group $G$ acting on $W$ with $\omega$ as only  fixed point, we have
an exact sequence
$$0\rightarrow \pi _1^{\omega}(W) \longrightarrow \pi _1^{\rm
o}(V)\longrightarrow G\rightarrow 0$$
(in particular $\pi _1^{\rm o}(V)=G$ if $W$ is smooth of dimension $\ge 2$).
Conversely, to each surjective homomorphism of $\pi _1^{\rm o}(V)$ onto a
finite group
$G$ corresponds an isolated singularity $(W,\omega)$ with an action of $G$
fixing only $\omega$ such that $W/G\cong V$; if $(V, {\rm o})$ is a
symplectic singularity, so is $(W,\omega)$. Therefore a first step in a
possible classification is to study isolated symplectic singularities
 with trivial local fundamental group. It turns out that the singularities
$(\overline{\cal O}_{\rm min},0)$ are of this type (with one exception):
\th Proposition
\enonce Let ${\goth g}$ be a simple complex Lie algebra, and ${\cal O}_{\rm
min}\i {\goth g}$ its minimal nilpotent orbit. Then 
$\pi _1^0(\overline{\cal O}_{\rm min})=0$ except if ${\goth g}$ is
 of type $C_r$ $(r\ge 1)$; in that case $\pi _1^0(\overline{\cal O}_{\rm
min})={\bf Z}/(2)$, and the corresponding double covering
of $\overline{\cal O}_{\rm min}$ is smooth.
\endth
\indp{\it Proof}: Consider the resolution $f:L^{-1}
\rightarrow\overline{\cal O}_{\rm min}$ (\ref{nil}); denote by
$E\i L^{-1} $ the zero section. Let $D$ be a tubular neighbourhood of $E$ in
$L^{-1} $, and $D^*=D\moins E$. Since the homogeneous space ${\bf P}{\cal
O}_{\rm min}$ is simply-connected, the homotopy exact sequence of the
fibration $f:D^*\rightarrow {\bf P}{\cal O}_{\rm min}$ reads
$$H_2({\bf P}{\cal O}_{\rm min},{\bf Z})\qfl{\partial}{\bf Z}\longrightarrow
\pi _1(D^*)\rightarrow 0\ ,$$where the map $\partial$ corresponds to the
Chern class $c_1(L^{-1} )\in H^2({\bf P}{\cal O}_{\rm min},{\bf Z})$.

\ind Put $\dim {\bf P}{\cal O}_{\rm min}=2r-1$. Since  $K_{{\bf P}{\cal
O}_{\rm min}}=L^{-r}$, the class
$c_1(L)$ is primitive unless ${\bf P}{\cal O}_{\rm min}={\bf P}^{2r-1}$, 
that is ${\goth g}$ is of type $C_r$. Assume this is not the case. The 
homotopy exact sequence gives
$\pi _1(D^*)=0$; since the pull back of any neighbourhood  of $0$ in
$\overline{\cal O}_{\rm min}$ contains a tubular neighbourhood of $E$, this
implies $\pi _1^0(\overline{\cal O}_{\rm min})=0$.
\ind If ${\goth g}$ is  of type
$C_r$ the same argument gives $\pi _1^0(\overline{\cal O}_{\rm min})={\bf
Z}/(2)$; actually $\overline{\cal O}_{\rm min}$ is the cone over
${\bf P}^{2r-1}$ embedded into ${\bf
P}^{(r+1)(2r-1)}$ by the Veronese embedding, and this cone is isomorphic to the
quotient of ${\bf C}^{2r}$ by the involution $v\mapsto -v$.\cqfd

\subsection  It would be interesting to find more examples of isolated
symplectic singularities with trivial local fundamental group, and also
examples with {\it infinite} local fundamental group.

\vskip2cm
\centerline{ REFERENCES} \vglue15pt\baselineskip12.8pt
\def\num#1{\smallskip \item{\hbox to\parindent{\enskip [#1]\hfill}}}
\parindent=1.3cm 
\num{B1} A.\ {\pc BEAUVILLE}: {\sl Fano contact manifolds and nilpotent 
orbits}.  Comment.\ Math.\ Helv.\ {\bf 73}, 566--583 (1998).
\num{B2} A.\ {\pc BEAUVILLE}: {\sl Riemannian Holonomy and Algebraic Geometry}.
Preprint math.AG/9902110.
\num{C-M} D.\ {\pc COLLINGWOOD}, W.\ {\pc MC}{\pc GOVERN}: {\sl Nilpotent 
orbits in semi-simple Lie algebras}. Van Nostrand Reinhold Co., New York
(1993).

\num{G} H.\ {\pc GRAUERT}: {\sl  \"Uber Modifikationen und exzeptionelle
analytische Mengen}.  Math.\ Ann.\ {\bf 146},  331--368 (1962). 
\num{K} Y.\ {\pc KIMURA}: {\sl On the hypersurfaces of Hermitian
symmetric spaces of compact type}. Osaka J.\ Math.\ {\bf 16},
97--119 (1979). 
 \num{L} C.\ {\pc LE}{\pc BRUN}: {\sl Fano manifolds, contact
structures, and quaternionic geometry}. Int.\ J.\ of Math. {\bf 6}, 
419--437 (1995).

\num{P} D.I.\ {\pc PANYUSHEV}: {\sl Rationality of singularities and the
Gorenstein properties of nilpotent orbits}. Functional\ Anal.\  Appl.\
{\bf 25}, 225--226 (1991). 
\num{R} M.\ {\pc REID}: {\sl Canonical $3$\tx folds}. Journ\'ees de G\'eometrie
Alg\'ebrique d'Angers (1979), 273--310; 
 Sijthoff \& Noordhoff (1980).
\num{W} J.\ {\pc WAHL}: {\sl  A cohomological characterization of}
${\bf P}^{n}$. Invent.\ Math.\ {\bf 72}, 315--322 (1983).

\vskip1cm
\def\pq#1{\eightrm#1\sixrm}
\hfill\vtop{\eightrm\hbox to 5cm{\hfill Arnaud {\pq BEAUVILLE}\hfill}
 \hbox to 5cm{\hfill DMI -- \'Ecole Normale
Sup\'erieure\hfill} \hbox to 5cm{\hfill (UMR 8553 du CNRS)\hfill}
\hbox to 5cm{\hfill  45 rue d'Ulm\hfill}
\hbox to 5cm{\hfill F-75230 {\pc PARIS} Cedex 05\hfill}}
\end